\documentclass[a4paper,11pt]{amsart}

\usepackage{cmap}            
\usepackage[T1]{fontenc}
\usepackage{lmodern}
\usepackage{amssymb}
\usepackage{mathtools}
\usepackage[
  pdfdisplaydoctitle=true,
  colorlinks=true,
  urlcolor=blue,
  citecolor=blue,
  linkcolor=blue,
  pdfstartview=FitH,
  pdfpagemode=UseNone,
  bookmarksnumbered=true
]{hyperref}

\newtheorem{thm}{Theorem}[section]

\newtheorem{lem}[thm]{Lemma}

\theoremstyle{definition}

\theoremstyle{remark}
\newtheorem{rem}[thm]{Remark}

\DeclarePairedDelimiter{\abs}{\lvert}{\rvert}
\DeclarePairedDelimiter{\norm}{\lVert}{\rVert}
\DeclarePairedDelimiter{\set}{\{}{\}}
\DeclarePairedDelimiter{\gen}{\langle}{\rangle}

\DeclareMathOperator{\GR}{GR}
\DeclareMathOperator{\Aut}{Aut}
\DeclareMathOperator{\card}{card}
\DeclareMathOperator{\interior}{int}
\DeclareMathOperator{\GL}{GL}
\DeclareMathOperator{\SL}{SL}

\newcommand{\Dp}[1]{\frac{\partial}{\partial #1}}
\newcommand{\C}{\mathbb{C}}
\newcommand{\Z}{\mathbb{Z}}
\newcommand{\M}{\mathbf{M}}
\newcommand{\Lt}{\Z[t,t^{-1}]}   

\hypersetup{
  pdftitle={Topological entropy and Burau representation},
  pdfauthor={Boris Kolev},
  pdfkeywords={Topological entropy, Burau representation, Fundamental group, Free differential calculus},
  pdfsubject={MSC: 20F36 (54C70, 57M05)},
  pdflang=en
}

\begin{document}

\title{Topological entropy and Burau representation}
\author[B. Kolev]{Boris Kolev}
\address{CNRS \& Université d'Aix-Marseille, CMI, 39 rue F. Joliot-Curie, 13453 Marseille cedex 13, France}
\email{kolev@cmi.univ-mrs.fr}


\date{August 25, 1989}
\commby{René Thom}

\begin{abstract}
  Let $f$ be an orientation-preserving homeomorphism of the disk $D^{2}$,
  $P$ a finite invariant subset and $[f_{P}]$ the isotopy class of $f$ in
  $D^{2}\setminus P$. We give a non-trivial lower bound of the topological
  entropy for maps in $[f_{P}]$, using the spectral radius of some
  specializations in $\GL(n,\C)$ of the Burau matrix associated with
  $[f_{P}]$, and we discuss some examples.
\end{abstract}

\maketitle


\section{Introduction}

Let $f$ be an orientation-preserving homeomorphism of the disk $D^{2}$ and
$P \subset \interior(D^{2})$ a finite invariant subset of $n$ points. Then
$f$ induces a homeomorphism $f_{P}$ of $D_{P} = D^{2}\setminus P$, and
$f_{P}$ induces an automorphism $f_{P\sharp}$ of
$\pi_{1}(D_{P}) \simeq F_{n}$, the free group on $n$ generators. If $h(f)$
is the \emph{topological entropy}~\cite{AKM1965} of $f$ and
$\GR(f_{P\sharp})$ is the \emph{growth rate} of $f_{P\sharp}$ on
$\pi_{1}(D_{P})$, Bowen~\cite{Bow1978} has shown that
\begin{equation}\label{BowenInequality}
  h(f) \geq \ln \GR(f_{P\sharp}).
\end{equation}
The problem is that, in general, $\GR(f_{P\sharp})$ is very hard to compute.

Let $f_{P*1}$ be the isomorphism induced by $f_{P}$ on the first homology
group $H_{1}(D_{P};\C)$. Then (see~\cite{Man1975})
\begin{equation}\label{ManningInequality}
  \GR(f_{P\sharp}) \geq R(f_{P*1}),
\end{equation}
where $R(f_{P*1})$ is the spectral radius of $f_{P*1}$. But $f_{P*1}$ is
just a permutation, so inequality~\eqref{ManningInequality} is always
trivial. Therefore, positive topological entropy cannot be detected this
way.

It is well known~\cite{Bir1974} that the subgroup of automorphisms of
$F_{n}$ induced by homeomorphisms of $D_{P}$ is isomorphic to
\emph{Artin's braid group}\footnote{To be perfectly precise, this result
  assumes that the base point of $\pi_{1}(D_{P})$ lies on the boundary of
  $D_{P}$.} $B_{n}$. There is a famous representation~\cite{Bir1974} of
$B_{n}$ in $\GL(n,\Lt)$, the \emph{Burau representation}. In this paper,
we prove the following result.

\begin{thm}\label{Main}
  Let $\alpha \in B_{n} \subset \Aut(F_{n})$ and let $B_{\alpha}$ be the
  Burau matrix of $\alpha$. Then
  \begin{equation}\label{KolevInequality}
    \GR(\alpha) \geq \sup \set*{R\bigl(B_{\alpha}(t)\bigr) \,;\, t\in \C,\ \abs{t} = 1},
  \end{equation}
  where $R(B_{\alpha}(t))$ is the spectral radius of the matrix
  $B_{\alpha}(t) \in \GL(n,\C)$ obtained by evaluating $B_{\alpha}$ at the
  complex number $t$.
\end{thm}

\begin{rem}
  This result is not really new, since it appears implicitly as a
  corollary of results of David Fried~\cite{Fri1986}. Nevertheless, we give
  here a direct proof, without any reference to ``twisted'' cohomology,
  using only the definition of the Burau representation and Fox's free
  differential calculus~\cite{CF1977}.
\end{rem}


\section{The braid group and the Burau representation}

\subsection{Artin's braid group}

There are several ways to introduce the braid group $B_{n}$. This
paragraph is only a short summary of~\cite{Bir1974}. We define $B_{n}$
algebraically as the subgroup of \emph{right automorphisms} $\alpha$ of the
free group
\begin{equation*}
  F_{n} = \gen{x_{1}, x_{2}, \dotsc, x_{n}}
\end{equation*}
which satisfy
\begin{align*}
  (x_{i})\alpha                  & = A_{i}\,x_{\mu_{i}}A_{i}^{-1}, \quad 1 \leq i \leq n, \\
  (x_{1}x_{2}\dotsm x_{n})\alpha & = x_{1}x_{2}\dotsm x_{n},
\end{align*}
where $(\mu_{1},\mu_{2},\dotsc,\mu_{n})$ is a permutation of
$(1,2,\dotsc,n)$ and $A_{i}$ is an element of $F_{n}$.

This group is finitely generated, with generators
\begin{equation}\label{generators}
  \sigma_{1}, \sigma_{2}, \dotsc, \sigma_{n-1},
\end{equation}
where
\begin{equation*}
  \left\{
  \begin{aligned}
    (x_{i})\sigma_{i}   & = x_{i}x_{i+1}x_{i}^{-1},       \\
    (x_{i+1})\sigma_{i} & = x_{i},                        \\
    (x_{j})\sigma_{i}   & = x_{j}, \quad j \neq i,\, i+1,
  \end{aligned}
  \right.
\end{equation*}
and relations
\begin{align}
  \sigma_{i}\sigma_{i+1}\sigma_{i} & = \sigma_{i+1}\sigma_{i}\sigma_{i+1}, \quad 1 \leq i \leq n-2, \label{relations} \\
  \sigma_{i}\sigma_{j}             & = \sigma_{j}\sigma_{i}, \quad \abs{i-j} \geq 2. \notag
\end{align}

\subsection{Free differential calculus}

Let $\Z F_{n}$ be the group ring of $F_{n}$ with integer coefficients. For
$j = 1,2,\dotsc,n$, there is a well-defined map
\begin{equation}\label{FreeDerivative}
  \Dp{x_{j}} \colon \Z F_{n} \to \Z F_{n},
\end{equation}
called the \emph{free derivative}, satisfying
\begin{align*}
  \Dp{x_{j}}(c_{1} + c_{2}) & = \Dp{x_{j}}(c_{1}) + \Dp{x_{j}}(c_{2}),      &  & c_{1}, c_{2} \in \Z F_{n}, \\
  \Dp{x_{j}}(g_{1} g_{2})   & = \Dp{x_{j}}(g_{1}) + g_{1}\Dp{x_{j}}(g_{2}), &  & g_{1}, g_{2} \in F_{n},    \\
  \Dp{x_{j}}(x_{i})         & = \delta_{ij},                                &  & 1 \leq i \leq n.
\end{align*}

Finally, let $\gen{t}$ be the free group generated by $t$, and let
\begin{equation}\label{projection}
  \varphi \colon F_{n} \to \gen{t}, \qquad x_{i} \mapsto t, \quad 1 \leq i \leq n.
\end{equation}
The morphism $\varphi$ extends to a ring morphism
$\Z F_{n} \to \Z\gen{t} = \Lt$, which we still denote by $\varphi$.

\subsection{Burau representation}

Consider an element $\alpha$ of $B_{n}$ and define the $n \times n$ matrix
$B_{\alpha} = (b_{ij})$ by
\begin{equation}\label{Burau}
  b_{ij} = \varphi\Bigl( \Dp{x_{j}}\bigl((x_{i})\alpha\bigr) \Bigr) \in \Lt,
  \qquad 1 \leq i,j \leq n.
\end{equation}
It is easy to show~\cite{CF1977} that
\begin{equation}\label{MatrixProduct}
  B_{\alpha\beta} = B_{\alpha}B_{\beta}.
\end{equation}
The map $\alpha \mapsto B_{\alpha}$ from $B_{n}$ to $\GL(n,\Lt)$ is the
\emph{Burau representation} of the braid group $B_{n}$.

Let $\M = \Lt^{n}$ and let $V_{1},V_{2},\dotsc,V_{n}$ denote the canonical
basis of the free module $\M$. A Burau matrix $B$ acts (on the right) on
row vectors of $\M$.

Any Burau matrix $B = (b_{ij})$ satisfies
\begin{equation}\label{Burau2}
  \sum_{k} b_{ik} = 1 \qquad (\forall i).
\end{equation}
Hence the submodule $\mathbf{H}$ of $\M$ defined by
\begin{equation}\label{Subspace}
  X_{1} + X_{2} + \dotsb + X_{n} = 0
\end{equation}
is invariant under $B$. The Burau representation is therefore reducible:
its restriction to $\mathbf{H}$ is an $(n-1)$-dimensional representation
$B^{r}$, called the \emph{reduced Burau representation}.

\begin{rem}
  There is an interesting relation between the Burau representation of a
  braid $\alpha$ and the Alexander polynomial, a famous invariant of links.
  If we close the braid $\alpha$ to obtain a link $\widetilde{\alpha}$ and
  let
  \begin{equation*}
    L = \widetilde{\alpha} \cup m,
  \end{equation*}
  where $m$ is the braid axis (a circle which encircles all the strings of
  the braid), then $L$ characterizes the link $\widetilde{\alpha}$ in the
  $3$-sphere $S^{3}$~\cite{Bir1974}, and the Alexander polynomial of $L$ is
  given by (cf.~\cite{Mor1985})
  \begin{equation}\label{Morton}
    \Delta(x,t) = \det\bigl(B_{\alpha}^{r}(t) - xI_{n-1}\bigr).
  \end{equation}
\end{rem}

\subsection{Spectral properties of Burau matrices}

A Burau matrix $B = (b_{ij})$ also satisfies
\begin{equation}\label{Burau1}
  \sum_{k} t^{k-1} b_{kj} = t^{j-1} \qquad (\forall j).
\end{equation}
Hence the row vector $(1,t,\dotsc,t^{n-1})$ is an eigenvector of $B$ with
eigenvalue $1$. If $P_{B}$ denotes the characteristic polynomial of $B$, we
thus have
\begin{equation}\label{Pol1}
  P_{B}(X) = \det(XI_{n} - B) = (X-1)\,P_{B^{r}}(X).
\end{equation}

Squier~\cite{Squ1984} has shown that the reduced Burau representation is
\emph{unitary}, in the sense that there exists a non-singular matrix
$J \in \GL(n-1,\Lt)$ such that
\begin{equation}\label{Squier}
  B^{r*}JB^{r} = J,
\end{equation}
where $B^{r*}$ is the transpose of the matrix obtained from $B^{r}$ by
exchanging $t$ and $t^{-1}$.

From equation~\eqref{Squier} we deduce that
\begin{equation}\label{Pol2}
  P_{B_{\alpha}^{r}}(1/X) = (-1/X)^{n-1}(-t)^{e}\,P_{B_{\alpha}^{r*}}(X),
\end{equation}
where $e$ is the algebraic sum of the exponents of the $\sigma_{i}$ in an
expression of $\alpha$.

Therefore, if $t \in \C$ with $\abs{t} = 1$, this relation becomes
\begin{equation}\label{Pol3}
  P_{B_{\alpha}^{r}}(1/X) = (-1/X)^{n-1}(-t)^{e}\,\overline{P_{B_{\alpha}^{r}}(\overline{X})}.
\end{equation}
Hence, if $\lambda$ is an eigenvalue of $B_{\alpha}^{r}$, then
$1/\overline{\lambda}$ is also an eigenvalue of $B_{\alpha}^{r}$, with the
same multiplicity.


\section{Growth rate of an automorphism of a finitely generated group}

Let $G$ be a finitely generated group, $S = \set{x_{1},x_{2},\dotsc,x_{n}}$
a generating set of $G$, and $L_{S}(g)$ the word length of $g \in G$ with
respect to $S$, that is, the minimal number of letters $x_{i}$ or
$x_{i}^{-1}$ needed to express $g$. The growth rate of an automorphism
$\alpha \colon G \to G$ is defined as
\begin{equation}\label{GrowthRate}
  \GR(\alpha) = \sup_{g\in G} \set*{\limsup_{p \to +\infty} \bigl( L_{S}((g)\alpha^{p})\bigr)^{1/p}}.
\end{equation}

It is not hard to see that $\GR(\alpha)$ does not depend on the generating
set $S$, that $1 \leq \GR(\alpha) < +\infty$, and that
\begin{equation}\label{GrowthRateNormal}
  \GR(\alpha) = \limsup_{p \to +\infty} \set*{\max_{1 \leq i \leq n} \bigl( L_{S}((x_{i})\alpha^{p})\bigr)^{1/p}}.
\end{equation}

Introduce the \emph{occurrence matrix} $A_{\alpha} = (a_{ij})$, where
$a_{ij}$ is the number of occurrences of $x_{j}^{\pm 1}$ in a reduced word
representing $(x_{i})\alpha$, and, for any complex matrix $M = (m_{ij})$,
the norm $\norm{M} = \max_{i} \sum_{j} \abs{m_{ij}}$. Since
$\norm{A_{\alpha^{p}}} = \max_{i} L_{S}((x_{i})\alpha^{p})$ and
$\norm{A_{\alpha^{p+q}}} \leq \norm{A_{\alpha^{p}}}\,\norm{A_{\alpha^{q}}}$,
equation~\eqref{GrowthRateNormal} gives (the limit exists by
submultiplicativity)
\begin{equation}\label{SpectralRadius}
  \GR(\alpha) = \lim_{p \to +\infty} \norm{A_{\alpha^{p}}}^{1/p}.
\end{equation}

\begin{rem}
  If $G = F_{n}$ is the free group, $S = \set{x_{1},x_{2},\dotsc,x_{n}}$ a
  system of free generators, $\alpha \in \Aut(F_{n})$, and if moreover
  there are no cancellations when we iterate $\alpha$, then
  $L_{S}((x_{i}x_{j})\alpha) = L_{S}((x_{i})\alpha) + L_{S}((x_{j})\alpha)$
  for $1 \leq i,j \leq n$. Hence $A_{\alpha^{p}} = A_{\alpha}^{p}$, and
  $\GR(\alpha)$ is just the spectral radius of the matrix $A_{\alpha}$. In
  the general case, we do not know any recurrence formula for
  $A_{\alpha^{p}}$.
\end{rem}


\section{Proof of the main theorem}

Let us come back to the case where $G = F_{n}$ and $\alpha \in B_{n}$.
Theorem~\ref{Main} is a consequence of the preceding considerations, as we
shall now see.

Let
$\omega = x_{\mu_{1}}^{\varepsilon_{1}}x_{\mu_{2}}^{\varepsilon_{2}} \dotsm x_{\mu_{r}}^{\varepsilon_{r}} \in F_{n}$
be a reduced word, with $\varepsilon_{k} = \pm 1$. We have
\begin{equation}\label{DerivativeOfOmega}
  \Dp{x_{j}}(\omega) = \sum_{k=1}^{r} \varepsilon_{k}\,\delta_{\mu_{k} j}\,
  x_{\mu_{1}}^{\varepsilon_{1}} \dotsm x_{\mu_{k-1}}^{\varepsilon_{k-1}}\,
  x_{\mu_{k}}^{(\varepsilon_{k}-1)/2}.
\end{equation}

Since $\omega$ is reduced, the group elements appearing in this sum are
pairwise distinct. Therefore, the number of occurrences of $x_{j}^{\pm 1}$
in $\omega$ is equal to the number of monomials in $\Dp{x_{j}}(\omega)$.
Moreover, if we write
\begin{equation*}
  \varphi\Bigl( \Dp{x_{j}}(\omega)\Bigr) = \sum_{m} c_{m}t^{m},
\end{equation*}
then the number of monomials in $\Dp{x_{j}}(\omega)$ is greater than or
equal to $\sum_{m} \abs{c_{m}}$. Now, taking $t \in \C$ with $\abs{t} = 1$,
we have
\begin{equation*}
  \abs[\Big]{\sum_{m} c_{m}t^{m}} \leq \sum_{m} \abs{c_{m}}.
\end{equation*}

Hence, if $\alpha \in B_{n}$, $A_{\alpha}$ is its occurrence matrix and
$B_{\alpha}$ its Burau matrix, we get, for all $t \in \C$ with
$\abs{t} = 1$,
\begin{equation*}
  a_{ij} \geq \abs{b_{ij}(t)}.
\end{equation*}

Applying this to $\alpha^{p}$ and using~\eqref{MatrixProduct}, we obtain
\begin{equation*}
  \norm{A_{\alpha^{p}}} \geq \norm{B_{\alpha^{p}}(t)} = \norm{B_{\alpha}(t)^{p}}.
\end{equation*}
Taking $p$-th roots and letting $p \to +\infty$, we conclude
with~\eqref{SpectralRadius} and Gelfand's formula
\begin{equation*}
  \lim_{p\to+\infty} \norm{M^{p}}^{1/p} = R(M),
\end{equation*}
which give $\GR(\alpha) \geq R(B_{\alpha}(t))$ for every $t$ with
$\abs{t} = 1$. This concludes the proof.

\begin{rem}
  For $t = 1$, $B_{\alpha}(1)$ is precisely the matrix of the map induced
  by $\alpha$ on the homology group $H_{1}(D_{P};\C)$, which is just a
  permutation matrix. Hence $t = 1$ always gives a trivial result.
\end{rem}


\section{Examples}

\subsection{\texorpdfstring{Example 1: $\alpha = \sigma_{1}\sigma_{2}^{-1} \in B_{3}$}{Example 1}}

The action on $F_{3}$ is given by
\begin{equation}\label{Example1}
  \left\{
  \begin{aligned}
    (x_{1})\alpha & = x_{1}x_{3}x_{1}^{-1}, \\
    (x_{2})\alpha & = x_{1},                \\
    (x_{3})\alpha & = x_{3}^{-1}x_{2}x_{3}.
  \end{aligned}
  \right.
\end{equation}

In this example, there are no cancellations when we iterate $\alpha$ on
each generator $x_{i}$. Hence, if $A_{\alpha^{p}}$ is the occurrence matrix
of $\alpha^{p}$, we have
\begin{equation*}
  A_{\alpha^{p}} = A_{\alpha}^{p}.
\end{equation*}
Therefore, we can compute $\GR(\alpha)$ explicitly as the spectral radius
of $A_{\alpha}$, and we find $(3 + \sqrt{5})/2$. On the other hand, the
characteristic polynomial of $B_{\alpha}$ is
\begin{equation*}
  P_{B_{\alpha}}(X) = (X-1)\bigl(X^{2} - (1-t-t^{-1})X + 1\bigr).
\end{equation*}
For $t = -1$, the root of largest modulus of $P_{B_{\alpha}}$ is precisely
$(3 + \sqrt{5})/2$. In this case, inequality~\eqref{KolevInequality} is an
equality.

This is, however, not surprising: $\GR(\alpha)$ is always equal to the
spectral radius of $B_{\alpha}(-1)$ when $\alpha \in B_{3}$. Let
$\overline{D}_{P}$ be the \emph{blow-up} of the $3$-punctured disk (where
each puncture is replaced by a circle~\cite{Bow1978}). The isotopy class
represented by $\alpha$ in $\overline{D}_{P}$ is periodic, reducible or
pseudo-Anosov (in the sense of Thurston~\cite{Thu1988}). In the first two
cases, $\GR(\alpha) = R(B_{\alpha}(-1)) = 1$. In the last case, the
invariant foliations of a pseudo-Anosov map $\phi$ in this class have
singularities only on the boundary of $\overline{D}_{P}$. The two-fold
cover of $\overline{D}_{P}$ associated with the hyperelliptic involution,
which we denote by $\overline{D}^{2}_{P}$, is a torus with $4$ holes (one
for each boundary component of $\overline{D}_{P}$). The lift
$\widetilde{\phi}$ of $\phi$ induces an Anosov diffeomorphism of the torus
obtained by collapsing each boundary curve to a point. The matrix
$B_{\alpha}^{r}(-1) \in \SL(2,\Z)$ is exactly the matrix of a linear Anosov
diffeomorphism $A$ in this class. Hence the entropy of $\phi$, which equals
that of $A$, is equal to $\ln R\bigl(B_{\alpha}(-1)\bigr)$.

\subsection{\texorpdfstring{Example 2: $\alpha = \sigma_{1}\sigma_{2}^{-1}\sigma_{3}^{-1} \in B_{4}$}{Example 2}}

The action on $F_{4}$ is given by
\begin{equation}\label{Example2}
  \left\{
  \begin{aligned}
    (x_{1})\alpha & = x_{1}x_{4}x_{1}^{-1}, \\
    (x_{2})\alpha & = x_{1},                \\
    (x_{3})\alpha & = x_{4}^{-1}x_{2}x_{4}, \\
    (x_{4})\alpha & = x_{4}^{-1}x_{3}x_{4}.
  \end{aligned}
  \right.
\end{equation}

In this example, there are cancellations when we iterate $\alpha$, so we
cannot compute $\GR(\alpha)$ directly. The characteristic polynomial of the
Burau matrix is
\begin{equation*}
  P_{B_{\alpha}}(X) = (X-1)\bigl(X^{3} - (1-t-t^{-1})X^{2} + (t^{-2}-t^{-1}+1)X + t^{-1}\bigr).
\end{equation*}

For $t = -1$, we get
\begin{equation*}
  P_{B^{r}_{\alpha}}(X) = (X-1)^{3},
\end{equation*}
which gives nothing. But for $t = j = e^{2\pi i/3}$ we obtain
\begin{equation*}
  P_{B^{r}_{\alpha}}(X) = X^{3} - 2X^{2} - 2j^{2}X + j^{2},
\end{equation*}
which has a root of modulus strictly greater than $1$. We can therefore
conclude that $h(f) > 0$ for all maps in the corresponding isotopy class.
In fact, this isotopy class is of pseudo-Anosov type. The growth rate
$\GR(\alpha)$ is the largest positive root of the polynomial
$X^{4} - 2X^{3} - 2X + 1$ (cf.~\cite[Appendix~A]{Kol1991}), namely
$\lambda = \bigl(1 + \sqrt{3} + \sqrt{2\sqrt{3}}\bigr)/2$, and we will show
that this number is not attained by our estimate~\eqref{KolevInequality}.

\begin{lem}
  In Example~2, we have
  \begin{equation*}
    \lambda > \sup_{\abs{t}=1} \set*{R\bigl(B_{\alpha}(t)\bigr)}.
  \end{equation*}
\end{lem}

\begin{proof}
  Throughout the proof, $t$ is a complex number of modulus $1$, and $P_{t}$
  denotes the cubic factor of the characteristic polynomial of
  $B_{\alpha}(t)$ computed above:
  \begin{equation*}
    P_{t}(X) = X^{3} - (1-t-t^{-1})X^{2} + (t^{-2}-t^{-1}+1)X + t^{-1}.
  \end{equation*}
  The eigenvalues of $B_{\alpha}(t)$ are thus the number $1$ and the three
  roots of $P_{t}$.

  We first show that no root of $P_{t}$ has modulus $\lambda$. Suppose, on
  the contrary, that $P_{t}(\lambda u) = 0$ for some complex number $u$ of
  modulus $1$. In other words, $u$ is a root of the polynomial
  \begin{equation*}
    F_{t}(X) = P_{t}(\lambda X) = \lambda^{3}X^{3} - (1-t-t^{-1})\lambda^{2}X^{2} + (t^{-2}-t^{-1}+1)\lambda X + t^{-1}.
  \end{equation*}
  Since $\overline{t} = t^{-1}$, $\overline{u} = u^{-1}$ and $\lambda$ is
  real, taking the complex conjugate of the equality $F_{t}(u) = 0$ and
  multiplying by $u^{3}$ shows that $u$ is also a root of
  \begin{equation*}
    G_{t}(X) = tX^{3} + (1-t+t^{2})\lambda X^{2} - (1-t-t^{-1})\lambda^{2}X + \lambda^{3}.
  \end{equation*}
  The polynomials $F_{t}$ and $G_{t}$ thus have a common root, so their
  resultant vanishes. Now this resultant can be computed explicitly.
  Reducing the powers of $\lambda$ by means of the relation
  $\lambda^{4} = 2\lambda^{3} + 2\lambda - 1$, one finds that it depends on
  $t$ only through the real number $z = t + t^{-1}$:
  \begin{gather*}
    \operatorname{Res}_{X}(F_{t},G_{t}) = c\,\Psi(z), \\
    \Psi(z) = 3z^{4} - 2(6-\sqrt{3})z^{3} - 12\sqrt{3}\,z^{2} + 6(4+3\sqrt{3})z + 4(12+7\sqrt{3}),
  \end{gather*}
  where $c = 2(834\lambda^{3} + 247\lambda^{2} + 568\lambda - 363)$ is a
  positive constant (note that $\sqrt{3} = \lambda + \lambda^{-1} - 1$).
  But $\Psi$ is positive on the real line. Indeed, expanding the right-hand
  side gives the identity
  \begin{gather*}
    3\,\Psi(z) = \bigl(3z^{2} - (6-\sqrt{3})z - 15\bigr)^{2} + 3\,q(z), \\
    q(z) = (17-8\sqrt{3})z^{2} + (28\sqrt{3}-36)z + (28\sqrt{3}-27),
  \end{gather*}
  and the quadratic polynomial $q$ takes only positive values: its leading coefficient is positive, because
  $17^{2} > 3\cdot 8^{2}$, and its discriminant $4(2043 - 1196\sqrt{3})$ is
  negative, because $2043^{2} < 3\cdot 1196^{2}$. The resultant is
  therefore non-zero, a contradiction.

  We can now conclude. By Theorem~\ref{Main}, and since
  $\GR(\alpha) = \lambda$, we have $R(B_{\alpha}(t)) \leq \lambda$. Equality
  would mean that some eigenvalue of $B_{\alpha}(t)$ has modulus $\lambda$.
  This eigenvalue cannot be $1$, because $\lambda > 1$, and we have just
  seen that it cannot be a root of $P_{t}$. Hence
  $R(B_{\alpha}(t)) < \lambda$ for every $t$ on the unit circle. Since the
  function $t \mapsto R(B_{\alpha}(t))$ is continuous and the unit circle
  is compact, its supremum is attained, and is therefore strictly less than
  $\lambda$.
\end{proof}

\subsection{\texorpdfstring{Example 3: $\alpha = \sigma_{4}\sigma_{3}\sigma_{2}\sigma_{1}\sigma_{4}\sigma_{3} \in B_{5}$}{Example 3}}

The action on $F_{5}$ is given by
\begin{equation}\label{Example3}
  \left\{
  \begin{aligned}
    (x_{1})\alpha & = x_{1}x_{2}x_{1}^{-1},                \\
    (x_{2})\alpha & = x_{1}x_{3}x_{4}x_{3}^{-1}x_{1}^{-1}, \\
    (x_{3})\alpha & = x_{1}x_{3}x_{5}x_{3}^{-1}x_{1}^{-1}, \\
    (x_{4})\alpha & = x_{1}x_{3}x_{1}^{-1},                \\
    (x_{5})\alpha & = x_{1}.
  \end{aligned}
  \right.
\end{equation}

This is another interesting example. Here again, $\alpha$ represents a
pseudo-Anosov class. The growth rate $\GR(\alpha)$ is the largest positive
root of $X^{4} - X^{3} - X^{2} - X + 1$, which we denote by $\lambda$
(cf.~\cite[Appendix~A]{Kol1991}).

On the other hand, we have
\begin{equation*}
  P_{B_{\alpha}(-1)}(X) = (X-1)\bigl(X^{4} + X^{3} - X^{2} + X + 1\bigr),
\end{equation*}
which has $-\lambda$ as a root. Hence, in this case,
\eqref{KolevInequality} is an equality. As in Example~1, the invariant
foliations of a pseudo-Anosov representative $\phi$ in this class have
singularities only on the boundary of $\overline{D}_{P}$. Since
$\card(P) = 5$ is odd, it can be shown that the lift $\widetilde{\phi}$ of
$\phi$ to the two-fold cover $\overline{D}_{P}^{2}$ has singularities only
on the boundary, and that they are of even order~\cite{Kol1991}. Hence the
invariant foliations of $\widetilde{\phi}$ are transversally orientable,
and $h(\widetilde{\phi}) = \ln R(\widetilde{\phi}_{*1})$. But
$B_{\alpha}^{r}(-1)$ is the homology matrix of $\widetilde{\phi}$, so
\begin{equation*}
  \ln \lambda = h(\widetilde{\phi}) = \ln R(\widetilde{\phi}_{*1}) = \ln R\bigl(B_{\alpha}^{r}(-1)\bigr).
\end{equation*}
This case of equality occurs whenever $\card(P)$ is odd and $\alpha$
represents a pseudo-Anosov class with singularities only on the boundary
of $\overline{D}_{P}$.

\begin{rem}
  Since this article was written, the subject has evolved and new results
  have been obtained (see for example~\cite{FL1999,LKS2002,TF2006,BB2007,RT2024}).
\end{rem}


\end{document}